\documentclass[final,a4paper]{article}
\usepackage[latin1]{inputenc}
\usepackage{fullpage}
\usepackage{amsmath}
\usepackage{xspace,url}
\usepackage[pdftex,colorlinks=true]{hyperref}
\usepackage[numbers]{natbib}
\renewcommand{\v}[1]{\boldsymbol{#1}} % Vectors and bivectors.
\renewcommand{\i}{\ensuremath{\v{i}}\xspace}
\renewcommand{\j}{\ensuremath{\v{j}}\xspace}
\renewcommand{\k}{\ensuremath{\v{k}}\xspace}
\newcommand{\Matrix}[1]{\ensuremath{\mathsf{#1}}}
\title{Canonic form of linear quaternion functions}
\author{Stephen J. Sangwine\thanks{Department of Computing and Electronic Systems,
                                   University of Essex, Wivenhoe Park, Colchester, CO4 3SQ,
                                   United Kingdom.
                                   Email:~\href{mailto:s.sangwine@ieee.org}{S.Sangwine@IEEE.org}
                                  }
                           \thanks{This work was partly funded by grant number EP/E010334/1 from the 
                                   United Kingdom Engineering and Physical Sciences Research Council.
                                  }
       }
\begin{document}
\maketitle
\begin{abstract}
The general linear quaternion function of degree one is a sum of terms with quaternion
coefficients on the left and right. The paper considers the canonic form of such a
function, and builds on the recent work of Todd Ell, who has shown that any such function
may be represented using at most four quaternion coefficients. In this paper, a new and
simple method is presented for obtaining these coefficients numerically using a matrix
approach which also gives an alternative proof of the canonic forms.
\end{abstract}
\section{Introduction}
In a recent paper \cite{arXiv:math/0702084v1}, Todd Ell has shown that a linear quaternion
function of the form:
\begin{equation}
\label{gqle}
f(q) = \sum_{p=1}^{P}m_p q\,n_p
\end{equation}
where the coefficients $m_p$ and $n_p$, and the variable $q$, are quaternions,
can be expressed in a canonic form using a maximum of four quaternion coefficients, irrespective
of the number of terms $P$. Ell's paper gives a method to obtain the four coefficients of the
canonic form algebraically without explicit decomposition of the original quaternion coefficients
$m_p$ and $n_p$ into their four components. This paper does not improve on that result, but it
does present a simple method for obtaining the four coefficients of the canonic form numerically.
The canonic form is\footnote{It is also possible to obtain a similar
form with the quaternion coefficients on the right and the quaternion operators $\i$, $\j$ and
$\k$ on the left, as was shown in \cite{arXiv:math/0702084v1}.}:
\begin{equation}
\label{canonic}
f(q) = Aq + Bq\i + Cq\j + Dq\k
\end{equation}
Subsequent to the publication of this result by Todd Ell, a paper by Littlewood in 1931 was
found to contain this canonic form \cite{10.1112/plms/s2-32.1.115} although Littlewood gave no
justification for his statement that it was `general' (assumed to mean `canonic')
nor a method to obtain the four coefficients from the general case as in equation \ref{gqle}.
No earlier papers are currently known that consider this form. In a later paper with Richardson,
the following equation appeared \cite[page 333]{10.1112/plms/s2-35.1.325}:
\[
a(x) = \sum_{q=1}^{n} = a_q x e_q = \sum_{p,q=1}^{n} a_{pq} e_p x e_q 
\]
The context of this equation is non-commutative algebras in general. In the specific case of quaternions
$e_0=1$, $e_1 =\i$, $e_2 = \j$ and $e_3 = \k$, and the result given in equation \ref{mex} in the next
section is a special case of Littlewood and Richardson's expression.

\section{Preliminaries}
\label{prelim}
Consider a single term in equation \ref{gqle}: $m_p q\,n_p$ where we write $m_p = w_p + x_p\i + y_p\j + z_p\k$ and
$n_p = w_p' + x_p'\i + y_p'\j + z_p'\k$. The product may be expanded as follows:
\begin{equation}
\label{mex}
\begin{aligned}
m_p q\,n_p &= (w_p + x_p\i + y_p\j + z_p\k)q(w_p' + x_p'\i + y_p'\j + z_p'\k)\\
           &= \begin{array}{*{4}{r@{\:+\:}}}
              w_p w_p'\, q & w_p x_p'\, q\i & w_p y_p'\, q\j & w_p z_p'\, q\k\\
              x_p w_p'\i q & x_p x_p'\i q\i & x_p y_p'\i q\j & x_p z_p'\i q\k\\
              y_p w_p'\j q & y_p x_p'\j q\i & y_p y_p'\j q\j & y_p z_p'\j q\k\\
              z_p w_p'\k q & z_p x_p'\k q\i & z_p y_p'\k q\j & z_p z_p'\k q\k\\
              \end{array}
\end{aligned}
\end{equation}
We can immediately see that the full function in equation \ref{gqle} can be written as a sum of
sixteen summations over the products of components of the left and right quaternion coefficients:
\begin{equation}
\label{smex}
\renewcommand{\S}{\displaystyle\sum_{p=1}^{P}}
f(q) = \S m_p q\,n_p =
\begin{array}{*{4}{r@{\:+\:}}}
\left(\S w_p w_p'\right)\, q & \left(\S w_p x_p'\right)\, q\i & \left(\S w_p y_p'\right)\, q\j & \left(\S w_p z_p'\right)\, q\k\\
\left(\S x_p w_p'\right)\i q & \left(\S x_p x_p'\right)\i q\i & \left(\S x_p y_p'\right)\i q\j & \left(\S x_p z_p'\right)\i q\k\\
\left(\S y_p w_p'\right)\j q & \left(\S y_p x_p'\right)\j q\i & \left(\S y_p y_p'\right)\j q\j & \left(\S y_p z_p'\right)\j q\k\\
\left(\S z_p w_p'\right)\k q & \left(\S z_p x_p'\right)\k q\i & \left(\S z_p y_p'\right)\k q\j & \left(\S z_p z_p'\right)\k q\k\\
\end{array}
\end{equation}
This expansion shows immediately that at least two canonic forms exist, as already
shown in \cite{arXiv:math/0702084v1}.
Grouping the four summations within the same column, each column sums to give a
quaternion coefficient on the left of $q$:
\[
f(q) = Aq + Bq\i + Cq\j + Dq\k
\]
An alternative is to group the four summations within the same row, where each row
sums to give a quaternion coefficient on the right of $q$:
\[
f(q) = qA' + \i qB' + \j qC' + \k qD'
\]
Now, a significant question is whether this is all, or whether there are additional canonic forms derivable
from this expansion. More accurately, are there additional canonic forms with 4 quaternion
coefficients (since a canonic form with more than 4 coefficients would be sub-optimal)?
The next section recasts the problem in matrix form in order to show how this question
might be answered.

In case this is thought by readers to be a trivial issue, Ljudmila Meister \cite[§3.5]{Meister:1997}
proposed the following as a general linear form: $f(q) = Aq + qB + CqD$. This has the correct number
of quaternion coefficients (four, the same as the canonic form given by Littlewood in 1931
\cite{10.1112/plms/s2-32.1.115}) and it looks very plausible as a canonic form, since there can
be only two terms with a single quaternion coefficient, one on the left and one on the right
(any others could be trivially combined with these since $Aq + Eq = (A+E)q$), and one double-sided
term with two coefficients. In fact, using the matrix formulation presented in the next section,
it is argued that this form cannot be canonic. It is very hard to show this by an algebraic argument,
but the matrix formulation presented in the next section makes it relatively simple.

\section{Matrix and outer product formulation}
The coefficients in equation \ref{mex} can be expressed as the following matrix, which is
the \emph{outer product} \cite{PenguinDictMaths3e} of two vectors representing the left and
right quaternion coefficients $m_p$ and $n_p$. It follows from the fact that this matrix is
the outer product of two vectors that it must be of rank 1 \cite{PenguinDictMaths3e}.
\begin{equation}
\label{outerprod}
\begin{pmatrix}
w_p w_p' & w_p x_p' & w_p y_p' & w_p z_p'\\
x_p w_p' & x_p x_p' & x_p y_p' & x_p z_p'\\
y_p w_p' & y_p x_p' & y_p y_p' & y_p z_p'\\
z_p w_p' & z_p x_p' & z_p y_p' & z_p z_p'\\
\end{pmatrix}
=
\begin{pmatrix}w_p \\ x_p \\ y_p \\ z_p\end{pmatrix}
\begin{pmatrix}w_p' & x_p' & y_p' & z_p\end{pmatrix}
\end{equation}
Similarly, we can express the coefficients in equation \ref{smex} as the following matrix.
It follows from the fact that this matrix is $4×4$ that it must be of rank 4 or less.
\begin{equation}
\renewcommand{\S}{\displaystyle\sum_{p=1}^{P}}
\Matrix{M} =
\begin{pmatrix}
\S w_p w_p' & \S w_p x_p' & \S w_p y_p' & \S w_p z_p'\\
\S x_p w_p' & \S x_p x_p' & \S x_p y_p' & \S x_p z_p'\\
\S y_p w_p' & \S y_p x_p' & \S y_p y_p' & \S y_p z_p'\\
\S z_p w_p' & \S z_p x_p' & \S z_p y_p' & \S z_p z_p'\\
\end{pmatrix}
=
\begin{pmatrix}
m_{11} & m_{12} & m_{13} & m_{14}\\
m_{21} & m_{22} & m_{23} & m_{24}\\
m_{31} & m_{32} & m_{33} & m_{34}\\
m_{41} & m_{42} & m_{43} & m_{44}\\
\end{pmatrix}
\end{equation}
\Matrix{M} may be decomposed into the sum of four matrices of rank 1, each of which can
be factored into an outer product of two vectors as in equation \ref{outerprod}.
It follows that, for any number
of terms $P$, equation \ref{gqle} can be expressed as \emph{at most} the sum of four terms,
that is $P=4$. The four rank 1 matrices which are summed to make \Matrix{M} can be obtained
from a singular value decomposition \cite{Golub:1996}, which can also yield the vectors
whose outer product gives one of the components of \Matrix{M}.
(Products of other pairs of columns will be zero because of orthogonality.)
This can be done as follows:
\[
\Matrix{M} = \Matrix{U}\Matrix{\Sigma}\Matrix{V}^{T} = \Matrix{L}\Matrix{V}^{T}
\]
where \Matrix{U} and \Matrix{V} are orthogonal, and \Matrix{\Sigma} is diagonal.
Corresponding columns of matrices \Matrix{L} and \Matrix{V} are vectors
whose outer products give the four matrices which sum to \Matrix{M}. These vectors are not unique,
although the singular values in \Matrix{\Sigma} are.

Clearly, a general linear quaternion function of the form in equation \ref{gqle} with an arbitrary
number of terms greater than 4 will yield a matrix \Matrix{M} of rank 4. It follows therefore that
any such function can be expressed using at most 4 terms, where the left and right quaternion
coefficients in each term are obtained from the outer product \emph{factorisation} of rank 1
matrices composing \Matrix{M} (these rank 1 matrices can be obtained from a singular-value
decomposition of \Matrix{M}):
\[
f(q) = AqE + BqF + CqG + DqH
\]
However, this is clearly not optimal, since it has eight quaternion coefficients, double the number
needed in the canonic form given in equation \ref{canonic}. However, by decomposing matrix \Matrix{M}
into the sum of four matrices each containing one column of \Matrix{M}, like this:
\begin{equation}
\renewcommand{\S}{\displaystyle\sum_{p=1}^{P}}
\Matrix{M} =
\begin{pmatrix}
m_{11} & 0 & 0 & 0\\
m_{21} & 0 & 0 & 0\\
m_{31} & 0 & 0 & 0\\
m_{41} & 0 & 0 & 0\\
\end{pmatrix}
+
\begin{pmatrix}
0 & m_{12} & 0 & 0\\
0 & m_{22} & 0 & 0\\
0 & m_{32} & 0 & 0\\
0 & m_{42} & 0 & 0\\
\end{pmatrix}
+
\begin{pmatrix}
0 & 0 & m_{31} & 0\\
0 & 0 & m_{32} & 0\\
0 & 0 & m_{33} & 0\\
0 & 0 & m_{34} & 0\\
\end{pmatrix}
+
\begin{pmatrix}
0 & 0 & 0 & m_{41}\\
0 & 0 & 0 & m_{42}\\
0 & 0 & 0 & m_{43}\\
0 & 0 & 0 & m_{44}\\
\end{pmatrix}
\end{equation}
we obtain the canonic form $f(q) = Aq + Bq\i + Cq\j + Dq\k$ in which the right-side
quaternion coefficients reduce to the degenerate set $1, \i, \j, \k$ and the left-side
coefficients are given by $A = m_{11} + m_{21}\i + m_{31}\j + m_{41}\k$ and so on.
(We could do the same with the rows of \Matrix{M}, in order to obtain the canonic form
with $1$, $\i$, $\j$ and $\k$ on the left.)

We now return to the `general' form proposed by Ljudmila Meister \cite[§3.5]{Meister:1997}
and mentioned in the last paragraph of §\ref{prelim}. This form is: $f(q) = Aq + qB + CqD$.
It should be clear from the foregoing that the three terms in this function would yield a
matrix \Matrix{M} of rank 3, and since the general case requires a matrix of rank 4, it is
not canonic. It might appear that adding a fourth term ($EqF$) would fix the problem, since
this would increase the rank to 4, but in fact this is not sufficient to yield a canonic
form. This can be seen as follows. $A$ corresponds to a column of the matrix \Matrix{M},
and $B$ to a row. This leaves three other columns and rows, which in the general case, will
be of rank 3. It follows that two double-sided terms like $CqD$ and $EqF$ are not sufficient to
represent the rest of the terms of a general linear function, since they cannot yield a
component of \Matrix{M} of rank greater than 2, and a rank of 3 is needed. We are thus led
to conclude that a function of the form $f(q) = Aq + qB + CqD + EqF + GqH$ is required,
which is again not optimal, since it has 8 quaternion coefficients. In fact, it is not
necessary for all the terms in this function to be full quaternions, but as we now show,
this does not improve on the form in equation \ref{canonic}.

If we take the first column and the first
row for $A$ and $B$, we would be left with a $3×3$ matrix, which in the general case
would be of rank 3. Factorising it into outer products we would need three terms with a
pure quaternion on each side. Noting also that it is not necessary for both of $A$ and $B$ to
be full quaternions, a `canonic' form appears to be:
$f(q) = Aq + qb + v_1 q v_2 + v_3 q v_4 + v_5 q v_6$, where $b$ and the $v_i$ are vectors
(pure quaternions). The number of real coefficients in this function is 25, 9 more than
in the canonic forms. However, when factorising the $3×3$ matrix within \Matrix{M} we can
factorise by columns, thus reducing the form to $f(q) = Aq + qb + v_1 q \i + v_3 q \j + v_5 q \k$
which has 16 real coefficients and is therefore canonic, although it is just a variation of
the form derived recently by Ell, and published by Littlewood in 1931.

It appears unlikely that there are any other decomposition(s) of \Matrix{M} that
would yield a canonic form with only 16 real coefficients.

\section{Conclusion}
In general, an arbitrary quaternion function of the form given in equation \ref{gqle} cannot
be reduced to less than four double-sided terms. That is, the function can always be reduced
to a sum of 4 double-sided terms, with 8 quaternion coefficients. There are useful special
cases where four of these coefficients are degenerate either on the left, or on the right
(but not mixed) as already shown in \cite{arXiv:math/0702084v1}.

The matrix formulation given in this paper provides a simple and systematic way to study the
problem further, although it appears unlikely that the canonic form in equation \ref{canonic}
can be improved on.

%\bibliographystyle{plainnat}
%\bibliography{sangwine,quaternions,maths,clifford,online}
\end{document}